\documentclass[12pt]{article}

\usepackage{amssymb}
\usepackage{amsthm}
\usepackage{authblk}
\usepackage{graphicx}

\pdfoutput=1

\newtheorem{thm1}{Theorem}

\newtheorem{lem1}{Lemma}

\newtheorem{fact1}{Fact}

\newtheorem{prop1}{Proposition}

\title{On the packing coloring of base-3 Sierpi\'{n}ski and $H$ graphs}

\author[1]{Fei Deng  }
\author[2]{Zehui Shao  }
\author[3]{Aleksander Vesel }
\affil[1]{School of Information Science and Technology,
Chengdu University of  Technology}
\affil[2]{Research Institute of Intelligence Software,Guangzhou University}
\affil[3]{Faculty of Natural Sciences and Mathematics, University of Maribor}

\begin{document}

\maketitle

\begin{abstract}
For a nondecreasing  sequence of integers $S=(s_1, s_2, \ldots)$ an $S$-packing $k$-coloring of a
graph $G$ is a mapping from $V(G)$ to $\{1, 2,\ldots,k\}$ such that vertices with color $i$ 
have pairwise distance greater than $s_i$.
By setting $s_i = d + \lfloor \frac{i-1}{n} \rfloor$ we obtain a  $(d,n)$-packing coloring of a graph $G$.
The smallest integer $k$ for which there exists a $(d,n)$-packing  coloring of $G$ is called the 
$(d,n)$-packing chromatic number of $G$. In the special case when $d$ and $n$ are both equal to one %$(d,n)=(1,1)$, 
we speak of the packing chromatic number of $G$. 
We determine the  packing chromatic number of the base-3 Sierpi\'{n}ski graphs $S_k$ and provide
new results on $(d,n)$-packing chromatic colorings, $d \le 2$, for this class of graphs.
By using a dynamic algorithm, we establish the packing chromatic number for  $H$-graphs $H(r)$. % if $r$ is odd and greater than or equal  to two.
%is and therefore establish the answer to the question
\end{abstract}

\section{Preliminaries}

A  {\em $k$-coloring} of a graph $G$ is a function $f$ from  $V(G)$ onto a set $ C = \{1, 2, \ldots, k \}$ (with no additional constraints). The elements of $C$ are called {\em colors}, while the set of vertices with the image (color) $i$ %imposed by $f$
is denoted by $X_i$.   %Let also $d_G(u, v)$ denote the usual shortest path distance between vertices $u$ and $v$ in a graph $G$. 
Let $u,v$ be vertices of a graph $G$. The distance between $u$ and $v$ in $G$, denoted by $d_G(u, v)$,
  equals the length of a shortest $u,v$-path (i.e. a path between $u$ and $v$)  
in $G$. %Note that $d_G(u, v)$  is also called the {\em distance between $u$ and $v$ in $G$}.

Let $f$ be a $k$-coloring of a graph $G$ with the corresponding sequence of color classes  $X_1, . . . , X_k$.
 If each color class $X_i$ is a set of vertices with the property that 
any distinct pair $u, v \in  X_i$ satisfies $d_G(u, v) > i$,
then $X_i$ is said to be an {\em $i$-packing}, while the sequence $X_1, . . . , X_k$ is called a {\em packing $i$-coloring}.
The smallest integer $k$ for which there exists a packing $k$-coloring of $G$
is called the {\em packing chromatic number of $G$} and it is denoted by  $\chi_{\rho}(G)$ \cite{BrKlRa, Goddard}.

A more general concept was  formally introduced in  \cite{GoddardXu2} as follows.
For a nondecreasing sequence
of integers $S = (s_1, s_2, \ldots)$, an {\em $S$-packing $k$-coloring} is a $k$-coloring $c$ of $V(G)$
such that for every $i$, with $1 \le i  \le k$, $c$ is an $s_i$-packing.
The {\em $S$-packing chromatic number}
of $G$ denoted by  $\chi_{\rho}^S(G)$, is the smallest $k$ such that $G$ admits an $S$-packing
$k$-coloring.

Gastineau et al. \cite{Gast} proposed the variation of the $S$-packing coloring, where for integers $n$ and $d$ the sequence $S = (s_1, s_2, \ldots)$  is given by $s_i = d + \lfloor \frac{i-1}{n} \rfloor$.
In this setting, an $S$-packing $k$-coloring of a graph $G$ is called a {\em $(d,n)$-packing $k$-coloring},
while  the smallest integer $k$ for which there exists a $(d,n)$-packing $k$-coloring of $G$
is called the {\em $(d,n)$-packing chromatic number} and denoted by
$\chi_{\rho}^{d,n}(G)$.

For $d=1$ and a sufficiently large $n$, a $(d,n)$-packing $k$-coloring is the classical graph coloring with $k$ colors.
A generalization of this observation gives the following
\begin{prop1} \label{chi}
Let $n \ge \ell$.
%Let $n$ be an integer.
%If $n \ge \chi(G) $, then $\chi_{\rho}^{1,n}(G) = \chi(G)$.
If $n \ge \chi_{\rho}^{d,\ell}(G) $, then $\chi_{\rho}^{d,n}(G) = \chi_{\rho}^{d,\ell}(G)$.
\end{prop1}

 Note also that a $(1,1)$-packing $k$-coloring is a packing $k$-coloring of a graph.

\section{Base-3 Sierpi\'{n}ski graphs}

Let $[k] := \{1,\ldots , k\}$ and $[k] _0:= \{0,\ldots , k-1\}$. The  base-3 Sierpi\'{n}ski graphs $S^k$ are defined such that
we start with $S^0 = K_1$. For $k \ge 1$, the vertex set of $S^k$ is $[3]_0^k$
and the edge set is defined recursively as

$E(S^k) = \{\{is, i t\} : i \in [3]_0, \{s, t\} \in E(S^{k-1})\} \cup
\{\{i j^{k-1}, ji^{n-1} \} | i, j \in [3]_0, i
= j \}$ .

\begin{figure}[hbt]
\centering

\unitlength 0.9mm % = 2.845pt
\linethickness{0.4pt}
\ifx\plotpoint\undefined\newsavebox{\plotpoint}\fi % GNUPLOT compatibility
\begin{picture}(96.072,80.249)(0,-5)
%\emline(49.403,83.249)(39.312,66.221)
\multiput(49.403,83.249)(-.0336358566,-.056760508){300}{\line(0,-1){.056760508}}
%\end
%\emline(23.966,39.732)(13.875,22.704)
\multiput(23.966,39.732)(-.0336358566,-.056760508){300}{\line(0,-1){.056760508}}
%\end
%\emline(73.999,39.522)(63.908,22.494)
\multiput(73.999,39.522)(-.0336358566,-.056760508){300}{\line(0,-1){.056760508}}
%\end
%\emline(36.789,61.385)(26.698,44.357)
\multiput(36.789,61.385)(-.0336358566,-.056760508){300}{\line(0,-1){.056760508}}
%\end
%\emline(11.352,17.869)(1.261,.841)
\multiput(11.352,17.869)(-.0336358566,-.056760508){300}{\line(0,-1){.056760508}}
%\end
%\emline(61.596,61.385)(51.505,44.357)
\multiput(61.596,61.385)(-.0336358566,-.056760508){300}{\line(0,-1){.056760508}}
%\end
%\emline(36.159,17.869)(26.068,.841)
\multiput(36.159,17.869)(-.0336358566,-.056760508){300}{\line(0,-1){.056760508}}
%\end
%\emline(86.192,17.659)(76.101,.631)
\multiput(86.192,17.659)(-.0336358566,-.056760508){300}{\line(0,-1){.056760508}}
%\end
\put(39.312,66.221){\line(1,0){19.971}}
\put(13.875,22.704){\line(1,0){19.971}}
\put(63.908,22.494){\line(1,0){19.971}}
\put(26.698,44.357){\line(1,0){19.971}}
\put(1.261,.841){\line(1,0){19.971}}
\put(51.295,.631){\line(1,0){19.971}}
\put(51.505,44.357){\line(1,0){19.971}}
\put(26.068,.841){\line(1,0){19.971}}
\put(76.101,.631){\line(1,0){19.971}}
%\emline(59.283,66.221)(49.403,83.249)
\multiput(59.283,66.221)(-.0337219552,.0581165611){293}{\line(0,1){.0581165611}}
%\end
%\emline(33.846,22.704)(23.966,39.732)
\multiput(33.846,22.704)(-.0337219552,.0581165611){293}{\line(0,1){.0581165611}}
%\end
%\emline(83.879,22.494)(73.999,39.522)
\multiput(83.879,22.494)(-.0337219552,.0581165611){293}{\line(0,1){.0581165611}}
%\end
%\emline(46.67,44.357)(36.789,61.385)
\multiput(46.67,44.357)(-.0337219552,.0581165611){293}{\line(0,1){.0581165611}}
%\end
%\emline(21.233,.841)(11.352,17.869)
\multiput(21.233,.841)(-.0337219552,.0581165611){293}{\line(0,1){.0581165611}}
%\end
%\emline(71.266,.631)(61.385,17.659)
\multiput(71.266,.631)(-.0337219552,.0581165611){293}{\line(0,1){.0581165611}}
%\end
%\emline(71.476,44.357)(61.596,61.385)
\multiput(71.476,44.357)(-.0337219552,.0581165611){293}{\line(0,1){.0581165611}}
%\end
%\emline(46.039,.841)(36.159,17.869)
\multiput(46.039,.841)(-.0337219552,.0581165611){293}{\line(0,1){.0581165611}}
%\end
%\emline(96.072,.631)(86.192,17.659)
\multiput(96.072,.631)(-.0337219552,.0581165611){293}{\line(0,1){.0581165611}}
%\end
\put(49.575,83.096){\circle*{2}}
\put(26.698,44.357){\circle*{2}}
\put(71.476,44.357){\circle*{2}}
\put(24.233,40.357){\circle*{2}}
\put(74.233,40.357){\circle*{2}}
\put(0.8,0.5){\circle*{2}}
\put(96.8,0.5){\circle*{2}}
\put(46,0.5){\circle*{2}}
\put(51,0.5){\circle*{2}}

\put(53.175,85.096){\makebox(0,0)[cc]{$0^{k}$}}
\put(23.698,48.357){\makebox(0,0)[cc]{$01^{k\!-\!1}$}}
\put(78,48.357){\makebox(0,0)[cc]{$02^{k\!-\!1}$}}

\put(20,42.357){\makebox(0,0)[cc]{$10^{k\!-\!1}$}}
\put(82,42.357){\makebox(0,0)[cc]{$20^{k\!-\!1}$}}

\put(42,-4){\makebox(0,0)[cc]{$12^{k\!-\!1}$}}
\put(55,-4){\makebox(0,0)[cc]{$21^{k\!-\!1}$}}

\put(0,5){\makebox(0,0)[cc]{$1^{k}$}}
\put(99,5){\makebox(0,0)[cc]{$2^{k}$}}

\put(49.192,70.686){\makebox(0,0)[cc]{{\small 00}$S^{k\!-\!2}$}}
\put(23.755,27.17){\makebox(0,0)[cc]{{\small 10}$S^{k\!-2\!}$}}
\put(73.789,26.96){\makebox(0,0)[cc]{{\small 20}$S^{k\!-\!2}$}}
\put(36.579,48.823){\makebox(0,0)[cc]{{\small 01}$S^{k\!-\!2}$}}
\put(11.142,5.307){\makebox(0,0)[cc]{{\small 11}$S^{k\!-\!2}$}}
\put(61.175,5.096){\makebox(0,0)[cc]{{\small 21}$\!S^{k\!-2\!}$}}
\put(61.385,48.823){\makebox(0,0)[cc]{{\small 02}$S^{k\!-\!2}$}}
\put(35.948,5.307){\makebox(0,0)[cc]{{\small 12}$S^{k\!-\!2}$}}
\put(85.982,5.096){\makebox(0,0)[cc]{{\small 22}$S^{k\!-\!2}$}}

\put(23.755,13.17){\makebox(0,0)[cc]{{\small 1}$S^{k-1}$}}
\put(73.789,13.17){\makebox(0,0)[cc]{{\small 2}$S^{k-1}$}}
\put(49.192,55){\makebox(0,0)[cc]{{\small 0}$S^{k-1}$}}

%\emline(39.312,66.01)(36.369,60.965)
\multiput(39.312,66.01)(-.03344474,-.05733385){88}{\line(0,-1){.05733385}}
%\end
%\emline(13.875,22.494)(10.932,17.449)
\multiput(13.875,22.494)(-.03344474,-.05733385){88}{\line(0,-1){.05733385}}
%\end
%\emline(63.908,22.284)(60.965,17.238)
\multiput(63.908,22.284)(-.03344474,-.05733385){88}{\line(0,-1){.05733385}}
%\end
%\emline(59.283,66.01)(61.596,61.175)
\multiput(59.283,66.01)(.03351399,-.0700747){69}{\line(0,-1){.0700747}}
%\end
%\emline(33.846,22.494)(36.159,17.659)
\multiput(33.846,22.494)(.03351399,-.0700747){69}{\line(0,-1){.0700747}}
%\end
%\emline(83.879,22.284)(86.192,17.449)
\multiput(83.879,22.284)(.03351399,-.0700747){69}{\line(0,-1){.0700747}}
%\end
\put(46.67,44.357){\line(1,0){5.676}}
\put(21.233,.841){\line(1,0){5.676}}
\put(71.266,.631){\line(1,0){5.676}}
%\emline(26.698,44.357)(23.966,39.732)
\multiput(26.698,44.357)(-.03373967,-.0570979){81}{\line(0,-1){.0570979}}
%\end
\put(51.295,.21){\line(0,1){0}}
\put(45.829,.631){\line(1,0){5.886}}
%\emline(60.965,17.238)(51.295,.631)
\multiput(60.965,17.238)(-.0336944557,-.0578665652){287}{\line(0,-1){.0578665652}}
%\end
%\emline(71.476,44.147)(74.209,39.312)
\multiput(71.476,44.147)(.03373967,-.05969326){81}{\line(0,-1){.05969326}}
%\end
\end{picture}
\caption{The structure of $S^k$}
\label{scierp}
\end{figure}
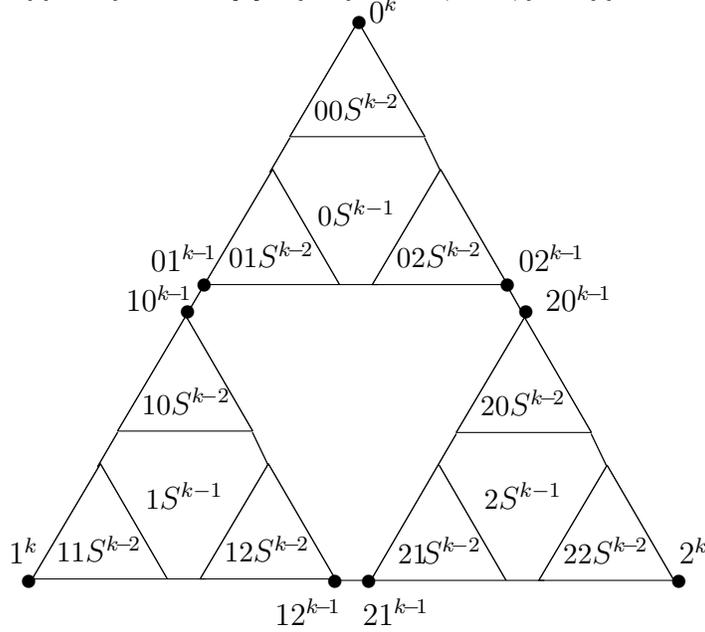

We can see that $S^k$ can be constructed from three copies of  $S^{k-1}$. More precisely, for every
$j \in [3]_0$ we make a copy of $S^{k-1}$ denoted by {\small \em  j}$S^{k-1}$ where
we concatenate $j$  to the left of each vertex in $S^{k-1}$. The construction is concluded by adding the edges: $\{0 1^{n-k}, 1 0^{k-1}\}$, $\{0 2^{k-1}, 2 0^{k-1}\}$, and $\{2 1^{k-1}, 1 2^{k-1}\}$ to the obtained graph. 
Obviously,   if $k \ge 2$, then $S^{k}$ is composed of nine copies of $S^{k-2}$ (see Fig. \ref{scierp}).

Note that if $u \in V(S^{k-1})$, then the corresponding vertex in  {\small \em  j}$S^{k-1}$ (a "copy" of $u$ in {\small \em  j}$S^{k-1}$ ) is of the form $ju$. 

Vertices of the form $i\ldots i=i^k$ are called {\em extreme  vertices}.
Clearly, if $k\ge 1$,  $S^k$ contains three extreme  vertices.

It is not difficult to establish the following (see also \cite{KZ})

\begin{fact1} \label{fact1}
Let $i,j \in [3]_0$ and $i \not = j$.
If $u \in V({\small  j}S^{k-1})$ and  $v \in V({\small  i}S^{k-1})$, then every shortest $u,v$-path
contains  vertices $ij^{k-1}$ and $ji^{k-1}$.
\end{fact1}

%The construction of $S^5$ is depicted in Fig. ??
%Note that in particular we showed that for $n \ge 6$ the graph $S^{n}$ is composed of $3^{n-5}$ copies of $S^{5}$.

Bre\v sar, Klav\v zar and Rall \cite{BKR}  showed the following
\begin{thm1} \label{sandi}
If $k \geq 5$, then $8 \le \chi_\rho(S^k) \le 9$.
\end{thm1}

In order to establish  the packing chromatic number of  base-3 Sierpi\'{n}ski graphs we need
the following definition.

\begin{figure}[hbt]
\centering
\includegraphics[width=4.2in]{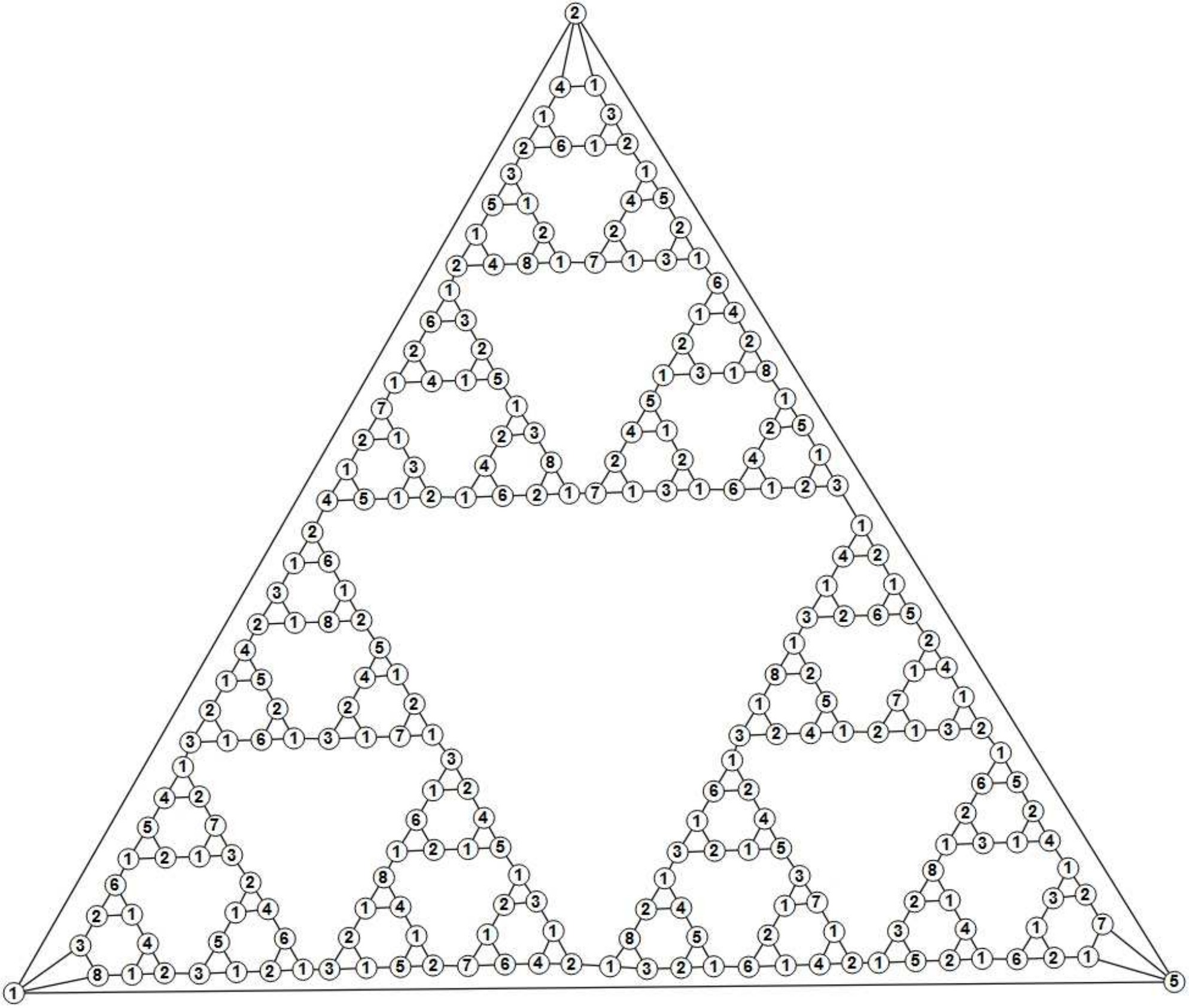}
\caption{a packing $8$-coloring of $T^5$}
\label{fig:T}
\end{figure}

Let $T^k$ be the graph obtained from $S^k$ by adding three edges that connect
its extreme vertices,
i.e. $V(T^k)=V(S^k)$ and $E(T^k)=E(S^k) \cup \{ \{0^k,1^k\}, \{2^k,1^k\},
\{0^k,2^k\} \}$. See for example Fig. \ref{fig:T} which shows $T^5$.

\begin{lem1} \label{Tgraph} Let $k \ge \ell$.
If $T^\ell$ admits an $(n,d)$-packing $b$-coloring such that
$d +\lfloor{ b -1 \over n } \rfloor \le 2^\ell$, then  $\chi_\rho^{d,n} (S^k) \le b$.
\end{lem1}

\begin{proof}
Let $f$ be an $(n,d)$-packing $b$-coloring of $T^\ell$.
Note that we show above that $S^k$ is composed of three copies of $S^{k-1}$.
For $k\ge \ell$ we  define  a $b$-coloring $f_k$  of $S^{k}$  as follows:

(i) $f_\ell(u) := f(u)$ for every $u \in V(S^\ell)$ (note that $V(T^\ell)=V(S^\ell)$.)

(ii) if $k  > \ell$,  then $f_k$ is obtained by applying $f_{k-1}$ to all three copies of $S^{k-1}$ in $S^{k}$.

We will show by induction on $k$ that $f_k$ is an $(n,d)$-packing $i$-coloring of $S^{k}$.

Since $S^\ell$ is a subgraph of $T^\ell$,
a $(n,d)$-packing $i$-coloring of $T^\ell$ is a $(n,d)$-packing $i$-coloring of $S^\ell$.
We therefore established that the claim holds for $k=\ell$.
 Let $u,v \in  V(S^{\ell+1})$. Thus, for
 some $x, y \in V(S^{\ell})$ and $i,j \in [3]_0$ we have $u=ix$ and $v=jy$.
The extreme vertices are connected in $T^\ell$, therefore, 
by Fact \ref{fact1}, we have $d_{S^{\ell+1}}(u,v)=d_{T^\ell}(x,y)$. Since $f(u)=f(x)$ and $f(v)=f(y)$, 
 the claim also holds for $k=\ell+1$. Let then $k\ge \ell+2$ and let assume that the claim holds for $S^{t}$, $t < k$.
 Let  $u,v \in V(S^k)$ such that $f_k(u)=f_k(v)$. We have to show that
 $d_{S^k}(u,v) > d + \lfloor {f_k(u) -1 \over n }\rfloor $.
  By the induction hypothesis, the restriction of $f_k$  to a copy
of $S^{k-1}$ (resp. a copy of $S^{k-2}$) is a $(n,d)$-packing $i$-coloring of the respective subgraph.
It follows that the claim clearly holds if $u$ and $v$ belong to the same copy of
$S^{k-1}$ (resp. $S^{k-2}$). If $u$ and $v$ belong to two copies
of $S^{k-2}$ which are not connected with an edge, then it is straightforward to see that
$d(u,v) > 2^{k-2} \ge 2^{\ell}$. Thus, this case is also settled.
 Finally, let $u$ and $v$
belong to two copies of $S^{k-2}$ which are connected with an edge.
We can say w.l.o.g. that for $x,y \in V(S^{k-2})$,
we have $u=01x$  and $v=10y$ or  $u=02x$  and $v=20y$ or  $u=12x$  and $v=21y$.

 If $u=01x$  and $v=10x$, we have $d_{S^k}(01x, 10y) = d_{S^k}(00x, 01y)$ and
$f_k(00x)=f_k(01y)=f_k(01x)$. By the inductive hypothesis,
$d_{S^k}(00x, 01y) > d +  \lfloor {f_k(00x) -1 \over n }\rfloor $  and this case is settled.
Since the proof for another cases is analogous,
we showed that $\chi_\rho^{d,n}(S^k) \le b$.
\end{proof}

\begin{thm1}

If $k \geq 5$, then $\chi_\rho(S^k) = 8$.
\end{thm1}

\begin{proof}
We found a packing $8$-coloring of $T^5$ depicted in  Fig. \ref{fig:T}. Thus,  by Lemma \ref{Tgraph}, we have
$\chi_\rho(S^k) = \chi_\rho^{1,1}(S^k) \le 8$. Since
Theorem \ref{sandi} yields the   lower bound, the proof is complete.
\end{proof}

Most of the other results of this section, $(d,n)$-packing chromatic numbers and bounds of  base-3 Sierpi\'{n}ski graphs,
 are  obtained by extensive computations which are based on ILP (proposed in \cite{zehvesel}) by Gurobi Optimizer 7.5 (http://www.gurobi.com/).
The results on $(d,n)$-packing chromatic numbers of $S^2$, $S^3$ and $S^4$ are given in Tables 1-3.

\begin{table}[hbt]
\caption{$(d,n)$-packing chromatic number of $S^2$} \label{S2packing}
% \begin{tabular}{p{1.2cm}p{0.8cm}p{0.6cm}p{1.2cm}p{0.8cm}p{0.6cm}p{1.2cm}p{0.6cm}p{0.8cm}p{1.2cm}p{0.8cm}p{0.3cm}}
 \begin{tabular}{p{1.cm}|p{1.1cm}p{1.1cm}p{1.2cm}p{1.1cm}p{1.1cm}p{1.2cm}}
\hline
$d \setminus n$ &  1 & 2 & 3 & 4 & 5 & 6  \\
\hline
1 & 5 & 3 & 3 & 3 & 3 & 3 \\
2 & 7 & 6 & 5  & 4 & 4 & 4 \\
$\geq 3$ & 9 & 9  & 9 & 9 & 9 & 9\\
\hline
 \end{tabular}
\end{table}

\begin{table}[hbt]
\caption{$(d,n)$-packing chromatic number of $S^3$} \label{S3packing}
% \begin{tabular}{p{1.2cm}p{0.8cm}p{0.6cm}p{1.2cm}p{0.8cm}p{0.6cm}p{1.2cm}p{0.6cm}p{0.8cm}p{1.2cm}p{0.8cm}p{0.3cm}}
 \begin{tabular}{p{1.cm}|p{0.6cm}p{0.6cm}p{0.6cm}p{0.6cm}p{0.6cm}p{0.6cm}p{0.6cm}p{0.6cm}p{0.6cm}p{0.6cm}p{0.6cm}p{0.6cm}}
\hline
$d \setminus n$ &  1 & 2 & 3 & 4 & 5 & 6 & 7  & 8 & 9 & 10& 11& 12   \\
\hline
1 & 7 & 4 & 3 & 3 & 3 & 3 & 3& 3& 3& 3& 3& 3\\
2 & 14 & 7  & 5  & 4  & 4 & 4 & 4 & 4 & 4& 4 & 4 & 4\\
3 & 19 & 13 & 10 &  9 &   9 & 9  & 9  & 9 & 9 & 9  & 9 & 9\\
4 & 21 & 17  & 14  &  11 &  11 &  11  &  11  &  10 &  10   &  10 &  10 &  10 \\
5 & 23 & 20 & 18  &  16 & 15 &  14 & 13 & 12 & 12& 12 & 12& 12\\
6 & 25 & 24 & 23 & 22 & 21 & 20 & 19 & 18 & 17& 16 & 16  & 16 \\
$\geq 7$ & 27 & 27& 27& 27& 27& 27& 27& 27& 27& 27& 27& 27\\
\hline
 \end{tabular}
\end{table}

\begin{table}[hbt]
\caption{$(d,n)$-packing chromatic number of $S^4$} \label{S4packing}
% \begin{tabular}{p{1.2cm}p{0.8cm}p{0.6cm}p{1.2cm}p{0.8cm}p{0.6cm}p{1.2cm}p{0.6cm}p{0.8cm}p{1.2cm}p{0.8cm}p{0.3cm}}
 \begin{tabular}{p{1.cm}|p{0.6cm}p{0.6cm}p{0.6cm}p{0.6cm}p{0.6cm}p{0.6cm}p{0.6cm}p{0.6cm}p{0.6cm}p{0.6cm}p{0.6cm}p{0.6cm}}
\hline
$d \setminus n$ &  1 & 2 & 3 & 4 & 5 & 6 & 7  & 8 & 9 & 10& 11& 12   \\
\hline
1 & 7 & 4 & 3 & 3 & 3 & 3 & 3& 3& 3& 3& 3& 3\\
2 & 20 & 7 & 6  & 4  & 4  & 4  & 4  & 4  & 4  & 4  & 4  & 4  \\
3 & 38 & 15  & 11 & 10 & 9 & 9 & 9 & 9 & 9 & 9 & 9 & 9 \\
4 & 46 & 20 & 15 & 12 & 12 & 11 &  11 & 11 & 11 & 10 & 10 & 10\\
5 & 54 & 30 & 20 & 18 &  16 & 14 & 14 & 13 & 13 & 13 & 13& 12 \\
6 & 60 & 41 & 27 & 23 & 22 & 21 & 20 & 20 & 19 & 18 & 18 & 17 \\
7 & 65 & 50 & 38 & 30 & 27 & 27 & 27 & 27 & 27 & 27 & 27 & 27    \\
8 & 67 & 54 & 44 &36 &31&29 &29 &29&28 &28&28&27 \\
9 & 69 & 58 & 50 & 43 & 37 & 33 & 32&32 &32 & 31& 31&31 \\
10 & 71 &62&56 &50&45&40&36&35&34&34&34&33\\
11 & 73 & 66 & 61 & 56 & 53 & 49 & 45 & 42 & 40 & 39 &37 &37\\
\hline
 \end{tabular}
\end{table}

\begin{thm1}
If $k \geq 2$ and $n \geq 2$, then
\begin{displaymath}
\chi_\rho^{1,n}(S^k)  =  \left\{ \begin{array}{ll}
 4,& \;k \geq 3 \;  and \; n=2 \\
 3,  &  \;  otherwise \\
 \end{array} \right.
\end{displaymath}
\end{thm1}

\begin{figure}[hbt]
\centering
\includegraphics[width=3.2in]{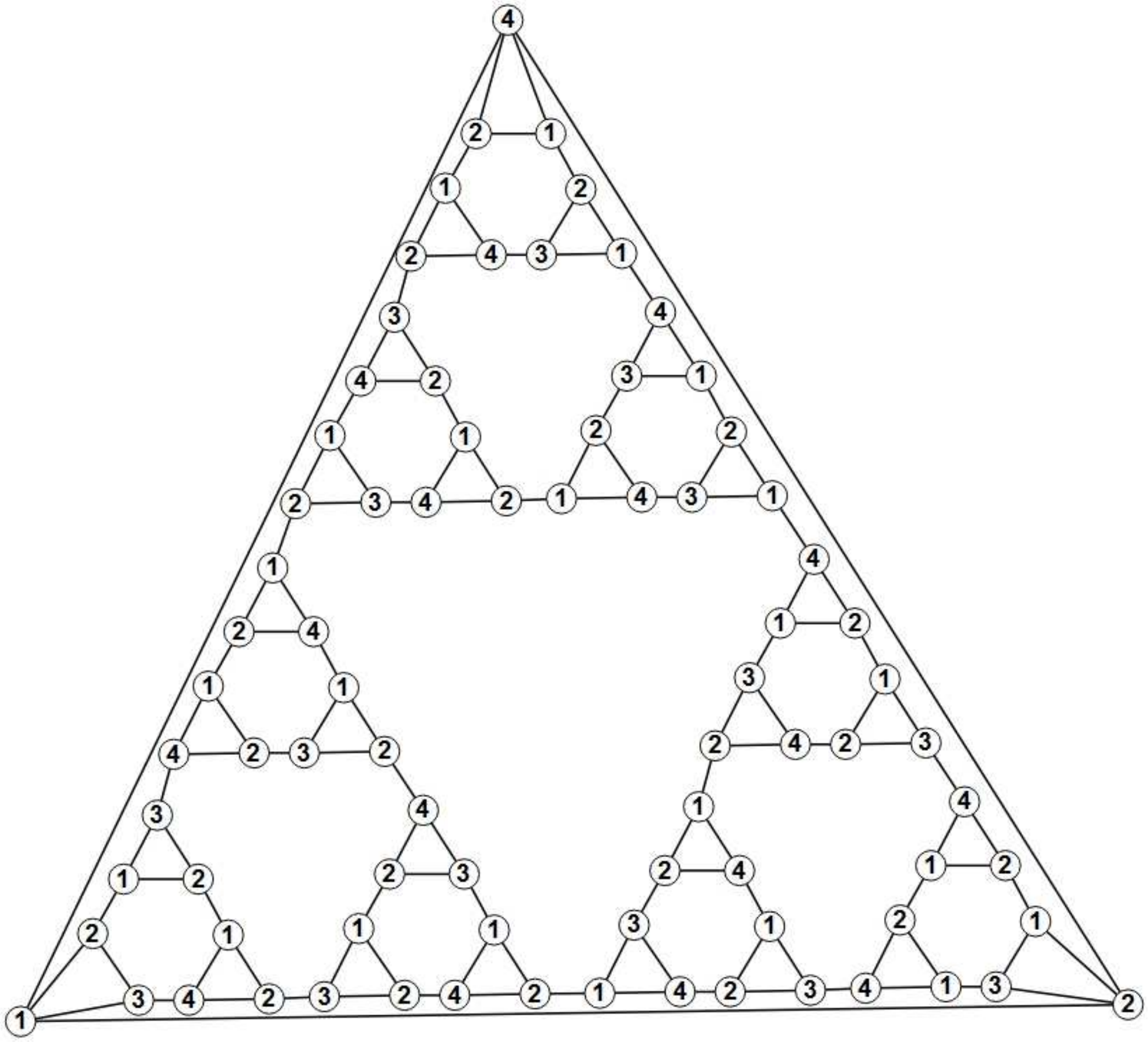}
\caption{A (1,2)-packing $4$-coloring of $T^4$}
\label{fig:T4}
\end{figure}

\begin{proof}
Since it is shown in \cite{JakKlav} that $\chi(S^k)=3$, from Proposition \ref{chi} it follows that
$\chi_\rho^{1,n}(S^k) =3$ for $n\ge 3$. Let then $n=2$.
For $k \in \{3,4,5 \}$, the values are established by a computer search.
Since $\chi_\rho^{1,2}(S^3) =4$, we have  $\chi_\rho^{1,2}(S^k) \ge 4$ for $k \ge 3$.
The upper bound is established by applying a (1,2)-packing $4$-coloring of $T^4$ depicted in Fig. \ref{fig:T4}
in Lemma \ref{Tgraph}.
\end{proof}

\begin{figure}[hbt]
\centering
\includegraphics[width=3.2in]{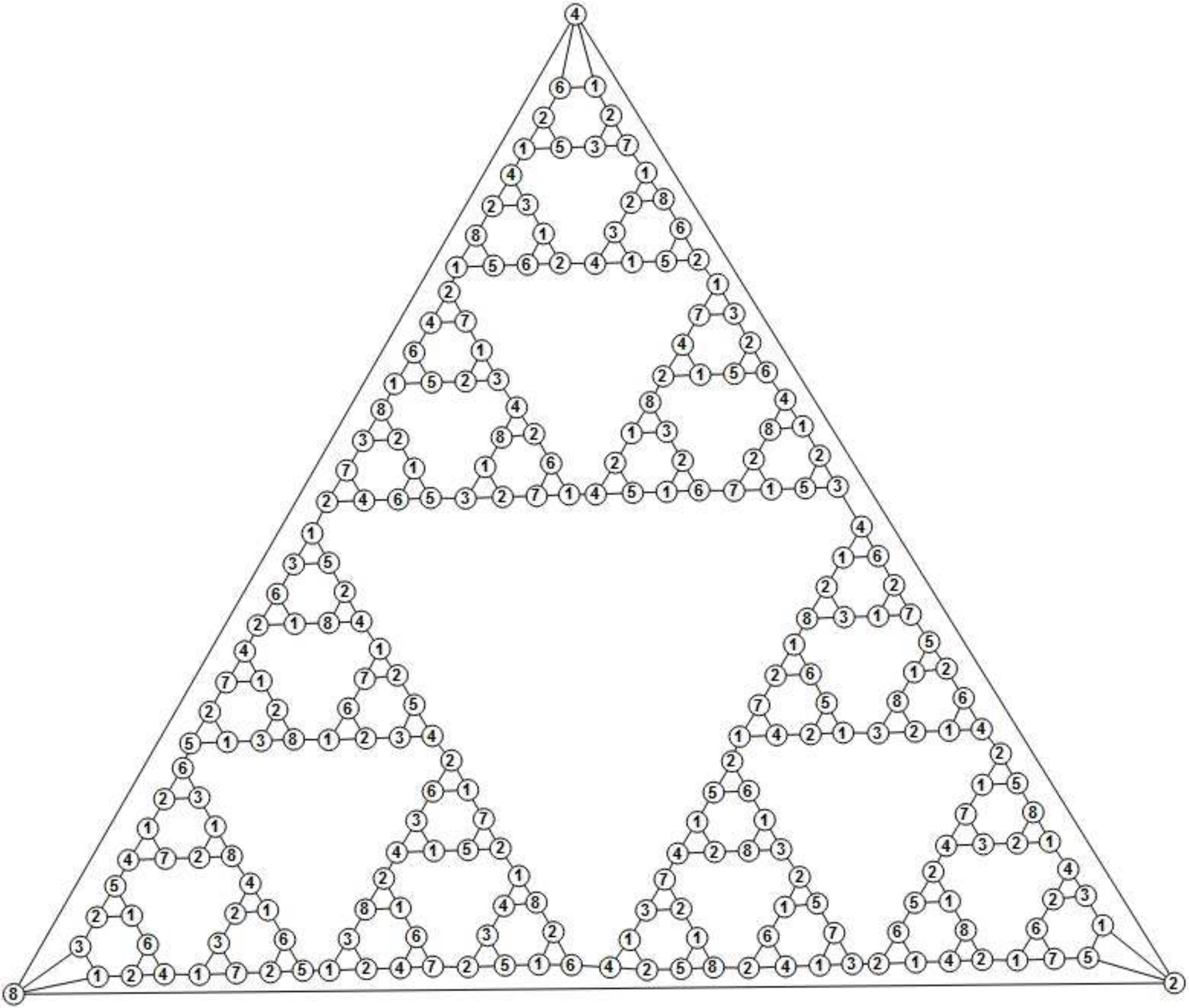}
\caption{A (2,2)-packing $8$-coloring of $T^5$}
\label{fig:T22}
\end{figure}

\begin{prop1}
%Let $k \geq 2$.
If $k \ge 6$ then $7 \le \chi_\rho^{2,2}(S^k) \le 8$.  Moreover,
\begin{displaymath}
\chi_\rho^{2,2}(S^k)  =  \left\{ \begin{array}{ll}
  6,  &  \;  k=2\\
  7,& \;k  \in \{3,4,5 \}\\
 % 7\; or\; 8,& \;k \geq 6\\
 \end{array} \right.
\end{displaymath}
\end{prop1}
%Moreover, if $k \ge 6$ then $7 \le \chi_\rho^{2,2}(S^k) \le 8$.  
\begin{proof}
For $k  \in \{2,3,4,5,6 \}$, the values are established by a computer search.
Since $\chi_\rho^{2,2}(S^3) =7$, we have  $\chi_\rho^{2,2}(S^k) \ge 7$ for $k \ge 3$.
The upper bound follows from a (2,2)-packing $8$-coloring of $T^5$ depicted in Fig. \ref{fig:T22}
and Lemma \ref{Tgraph}.
\end{proof}

%I guess the graph $S^6$ admits no (2,2)-packing $7$-coloring, since I have computed this case for about two days and find no such a coloring.

\begin{prop1}
%Let $n \geq 2$.
\begin{displaymath}
\chi_\rho^{2,3}(S^k)  =  \left\{ \begin{array}{ll}
 5,  &  \;  k  \in \{2,3 \} \\
 6,& \;k \geq 4\\
 \end{array} \right.
\end{displaymath}
\end{prop1}

\begin{figure}[hbt]
\centering
\includegraphics[width=3.2in]{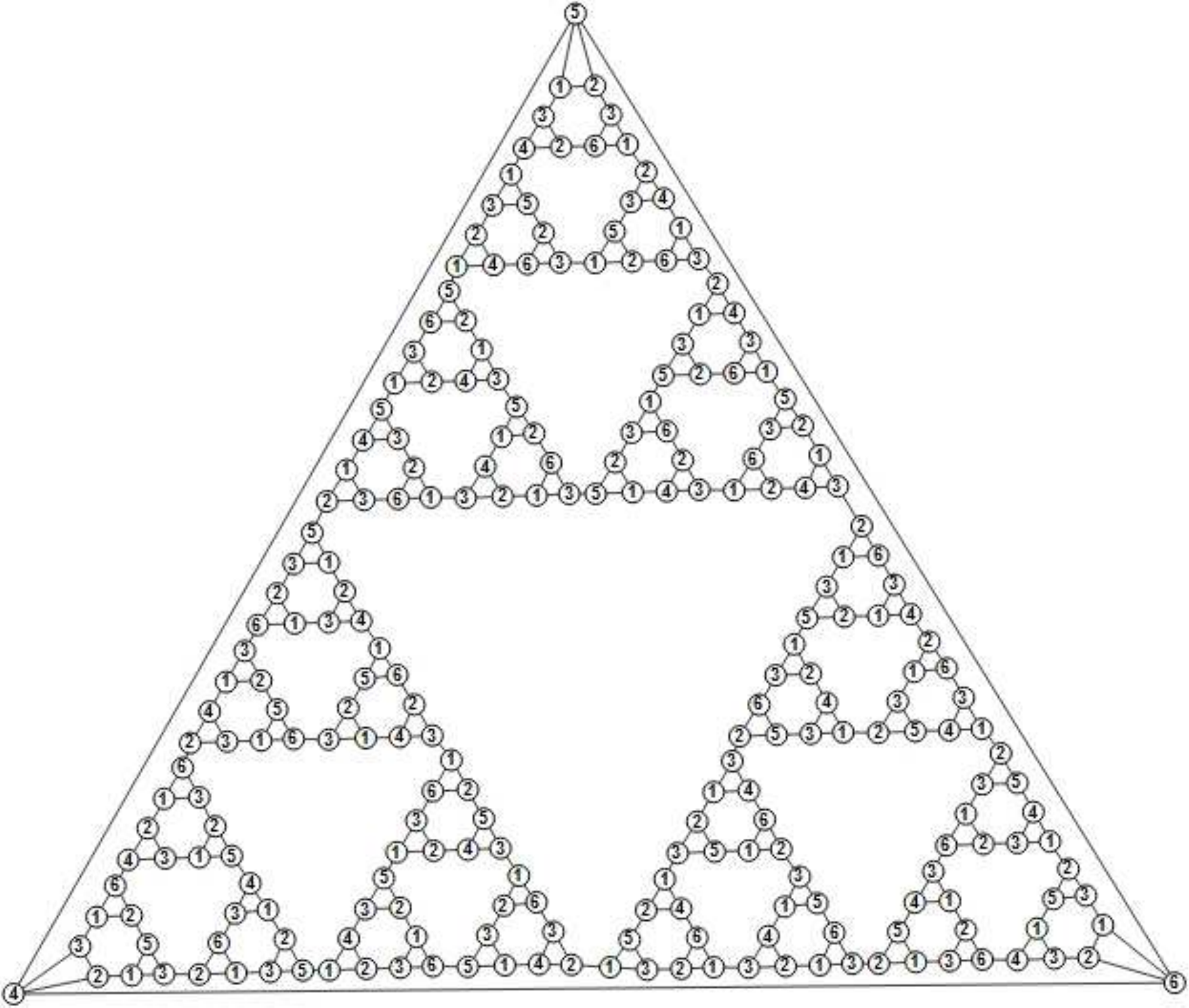}
\caption{A (2,3)-packing $6$-coloring of $T^5$}
\label{fig:T23}
\end{figure}

\begin{proof}
For $k  \in \{2,3 \}$, the values are established by a computer search.
Since $\chi_\rho^{2,3}(S^4) =6$, we have  $\chi_\rho^{2,3}(S^k) \ge 6$ for $k \ge 4$.
The upper bound is established by using a (2,3)-packing $6$-coloring of $T^5$ depicted in Fig. \ref{fig:T23}
in Lemma \ref{Tgraph}.
\end{proof}

\section{$H$-graphs}

The $H$-graph $H(r)$, $r \ge 2$, is the 3-regular graph of order $6r$, with vertex set
$V (H(r)) = \{ u_i,v_i,w_i : 0 \le i \le 2r - 1 \}$
and edge set $E(H(r)) = \{(u_i, u_{i+1}), (w_i, wi_{i+1}), (u_i, v_i), (v_i,w_i) : 0 \le i \le 2r -1 \} \cup
\{ (v_{2i}; v_{2i+1}) : 0 \le i \le r - 1 \}$, where subscripts are taken modulo $2r$.

For $i \in [r]_0$ the set of vertices $H_{r,i} = \{u_i, u_{i+1}, v_i, v_{i+1}, w_i, w_{i+1} \}$ induce a subgraph of $H(r)$ called
an {\em $i$-th column} of $H(r)$.

La\"{i}che and Sopena \cite{sopena} showed
\begin{thm1} \label{sop}
For every integer $r \ge 2$, $\chi_{\rho}(H(r)) = 5$   if r is even, and $6 \le \chi_{\rho}(H(r)) \le 7$  if $r$ is
odd.
\end{thm1}

The question  whether it is true that for $r \ge 3$, $r$ odd,  $H(r) \le 6$ for every $H$-graph $H(r)$ is posed
in \cite{sopena}.
In order to answer this question we use a dynamic algorithm which
was introduced  (in a very general framework) in \cite{klve-01}.
In this paper, this concept is applied for searching the  packing chromatic number of  $H$-graphs.  %$H(r)$, $r$ odd, as follows.

Consider the graph $H_P(r)$ obtained from $H(r)$ by removing the edges $u_0u_{2r-1}$ and $w_0w_{2r-1}$ 
(see for example $H_P(3)$ depicted in Fig. \ref{H}).
Then the vertices of the directed graph $D$ are all 6-packing colorings of  $H_P(3)$.
Let $f, g \in V(D)$, i.e. $f$ and $g$ are  6-packing colorings which assign a color to every vertex of $H_P(3)$.

If  $f, g \in V(D)$, then let $f_g$ be a 6-coloring  of  $H_P(4)$ such that
$f_g(H_{3,i}) :=  f(H_{3,i})$, $i \in [3]_0$, and  $f_g(H_{3,3}) :=  g(H_{3,2})$, i.e. $f_g$ is composed of $f$ and the last column of $g$.
We make an arc from $f$ to $g$ in $D$   if and only if the following two conditions are fullfiled:

(i) $f(H_{3,1}) = g(H_{3,0})$ and  $f(H_{3,2}) = g(H_{3,1})$, i.e., the coloring of the $1$st (resp.  $2$nd) column of $f$
coincide with the coloring of the $0$th (resp.  $1$st) column of $g$.

(ii) $f_g$ is a 6-packing coloring of  $H_P(4)$.

\begin{figure}[hbt]
\centering

\unitlength 0.8mm % = 2.85pt
\linethickness{0.4pt}
\ifx\plotpoint\undefined\newsavebox{\plotpoint}\fi % GNUPLOT compatibility
\begin{picture}(100,35.35)(9,0)
\put(8,5){\circle*{2.}}
\put(28,5){\circle*{2.}}
\put(48,5){\circle*{2.}}
\put(68,5){\circle*{2.}}
\put(88,5){\circle*{2.}}
\put(108,5){\circle*{2.}}

\put(8,20){\circle*{2.}}
\put(28,20){\circle*{2.}}
\put(48,20){\circle*{2.}}
\put(68,20){\circle*{2.}}
\put(88,20){\circle*{2.}}
\put(108,20){\circle*{2.}}

\put(8,35){\circle*{2.}}
\put(28,35){\circle*{2.}}
\put(48,35){\circle*{2.}}
\put(68,35){\circle*{2.}}
\put(88,35){\circle*{2.}}
\put(108,35){\circle*{2.}}

\put(5,2){\makebox(0,0)[cc]{$w_0$}}
\put(25,2){\makebox(0,0)[cc]{$w_1$}}
\put(45,2){\makebox(0,0)[cc]{$w_2$}}
\put(65,2){\makebox(0,0)[cc]{$w_3$}}
\put(85,2){\makebox(0,0)[cc]{$w_4$}}
\put(105,2){\makebox(0,0)[cc]{$w_5$}}

\put(5,17){\makebox(0,0)[cc]{$v_0$}}
\put(25,17){\makebox(0,0)[cc]{$v_1$}}
\put(45,17){\makebox(0,0)[cc]{$v_2$}}
\put(65,17){\makebox(0,0)[cc]{$v_3$}}
\put(85,17){\makebox(0,0)[cc]{$v_4$}}
\put(105,17){\makebox(0,0)[cc]{$v_5$}}

\put(5,32){\makebox(0,0)[cc]{$u_0$}}
\put(25,32){\makebox(0,0)[cc]{$u_1$}}
\put(45,32){\makebox(0,0)[cc]{$u_2$}}
\put(65,32){\makebox(0,0)[cc]{$u_3$}}
\put(85,32){\makebox(0,0)[cc]{$u_4$}}
\put(105,32){\makebox(0,0)[cc]{$u_5$}}

\put(8,5){\line(1,0){20}}
\put(28,5){\line(1,0){20}}
\put(48,5){\line(1,0){20}}
\put(68,5){\line(1,0){20}}
\put(88,5){\line(1,0){20}}

\put(8,20){\line(1,0){20}}
\put(48,20){\line(1,0){20}}
\put(88,20){\line(1,0){20}}

\put(8,35){\line(1,0){20}}
\put(28,35){\line(1,0){20}}
\put(48,35){\line(1,0){20}}
\put(68,35){\line(1,0){20}}
\put(88,35){\line(1,0){20}}

\put(8,5){\line(0,1){15}}
\put(28,5){\line(0,1){15}}
\put(48,5){\line(0,1){15}}
\put(68,5){\line(0,1){15}}
\put(88,5){\line(0,1){15}}
\put(108,5){\line(0,1){15}}

\put(8,20){\line(0,1){15}}
\put(28,20){\line(0,1){15}}
\put(48,20){\line(0,1){15}}
\put(68,20){\line(0,1){15}}
\put(88,20){\line(0,1){15}}
\put(108,20){\line(0,1){15}}

\end{picture}

\caption{ $H_P(3)$}
\label{H}
\end{figure}

\begin{lem1}
\label{basic}
Let $r \ge 4$ be an integer.  Then
$H(r)$ admits a  packing $6$-coloring if and only if $D$  contains a closed directed walk of length $r-2$.
\end{lem1}

\begin{proof}
Suppose first that $D$  contains  a closed directed walk $P$ of length $r-2$. Note that every arc of $D$ corresponds to a
 a 6-coloring  of  $H_P(4)$. Analogously, arcs of $P$ correspond to a 6-coloring  $\psi$ of  $H(r)$. We have to show $\psi$ is a packing 6-coloring  of  $H(r)$, i.e. that for every $u,v \in  H(r)$ with 
 $\psi(u)=\psi(v)$ we have $d_{H(r)}(u,v) > \psi(u)$. Let $u$ and $v$ belong to $i$-th and 
$j$-th column of $H(r)$, $i\ge j$, respectively. If $i-j \ge 4$, then  $d_{H(r)}(u,v) > 6$ and we are done. If $i-j \le 3$, then the colorings of $H_{r,j},H_{r,j+1},H_{r,j+2},H_{r,j+3}$ correspond to 
an arc of $D$. Thus, $d_{H(r)}(u,v) > \psi(u)$.
 \iffalse
 Let $f, g$ and $h$ be consecutive vertices of $P$.
Note that $f_{g}$ and $g_{h}$ are both 6-packing colorings of  $H_P(4)$.
Let $f_{gh}$ be a 6-coloring of  $H_P(5)$ such that
$f_{gh}(H_{3,i}) :=  f(H_{3,i})$, $i \in [3]_0$, and  $f_{gh}(H_{3,j}) :=  h(H_{3,j})$, $j \in \{ 1,2\}$,
i.e. $f_{gh}$ is composed of  $f$ together with the last two columns of $h$.
Let $x$ and $y$ be vertices of a copy of $H_P(5)$ that corresponds to $f_{gh}$ such that $f(x)=f(y)$.
We may assume that $x$ and $y$ belong to the $i$-th and the $j$-th column of $H_P(5)$, $i,j \in [5]_0$, $i \le j$, respectively.
We have to show that $d(x,y) > f(x)$. If $i=0$ and $j=5$, then   $d(x,y) \ge 7$ and we are done.
Otherwise, $x$ and $y$ both belong to  either to $[4]_0$ or  to $\{ 1,2,3,4\}$. 
But $f_{g}$ and $g_{h}$ are both 6-packing colorings of  $H_P(4)$ and the claim clearly holds.
\fi
This observation concludes the first part of the proof.

Assume now that $H(r)$ admits a  packing $6$-coloring denoted by $\psi$.
By the definition of $D$, the restriction of  $\psi$  to a copy of $H_P(3)$ (resp. a copy of $H_P(4)$) in $H(r)$ corresponds to a vertex (resp. an arc) of $D$. It follows that $\psi$ corresponds to a closed directed walk of length $r-2$.
\end{proof}

\begin{thm1}
If $r \geq 2$, then
\begin{displaymath}
\chi_{\rho}(H(r))   =  \left\{ \begin{array}{ll}
 5,& r \; is \; even  \\
 7,  &  \;  otherwise \\
 \end{array} \right.
\end{displaymath}
\end{thm1}

\begin{proof}
By Theorem \ref{sop}, we have to show that $\chi_{\rho}(H(r)) \ge 7$ if $r$ is odd.
In order to prove  this, we compute the directed graph $D$ (with 8336 vertices) as described above.
Since the computer program find no closed walk of odd length, by Lemma \ref{basic}  the assertion follows.
\end{proof}

%It is not difficult to compute $\chi_\rho^{1,n}(S^k)$ for every $k$ and $n$ (see the attachment).
%It is maybe also possible to compute $\chi_\rho^{2,n}(S^k)$ for every $k$ and $n$.
%It seems that $\chi_\rho^{2,n}(S^k) = 4$ if $n \geq 4$.
%This result could maybe be confirmed with a (2,4)-packing 4-coloring for some $T^k$.

\section*
{Acknowledgement}
%{\it Acknowledgement.}
This work is supported by the National Key Research and
Development Program
under grants 2017YFB0802300 and 2017YFB0802303,
the National Natural Science Foundation of China under the grant 11361008, 
the Applied Basic Research
(Key Project) of Sichuan Province under grant 2017JY0095 and 
the Ministry of Science of Slovenia under the grants P1-0297 and J1-7110.

\end{document}